\documentclass[12pt,reqno]{article}
\usepackage{graphicx}
\usepackage{amsmath}
\usepackage{amsthm}
\usepackage{latexsym,bm}
\usepackage{amssymb}
\usepackage{indentfirst}
\newtheorem{thm}{Theorem}[section]

\newtheorem{lem}[thm]{Lemma}

\numberwithin{equation}{section}
\begin{document}

\title{{\bf  A Riemannian Bieberbach estimate
}}

\author{{\bf F. Fontenele\thanks{Work partially supported by CNPq (Brazil)} \  and  F. Xavier}}

\date{}
\maketitle

\begin{quote}
\small {\bf Abstract}. The Bieberbach estimate, a pivotal result
in the classical theory of univalent functions, states that any
injective holomorphic function $f$ on the open unit disc $D$
satisfies $|f''(0)|\leq 4 |f'(0)|$.  We  generalize the Bieberbach
estimate by proving a version of the inequality that applies to
all injective smooth conformal immersions $f : D\to \Bbb R^n,\,
n\geq 2$. The new estimate involves two correction terms. The
first one is geometric, coming from the second fundamental form of
the image surface $f(D)$. The second term
is of a dynamical nature,  and involves
certain Riemannian quantities associated  to conformal attractors.
Our results are partly motivated by a conjecture in the theory of  embedded minimal surfaces.%The results in this paper are partly motivated by the following
%geometric conjecture: \it Any embedded minimal surface in $\Bbb
%R^3$ that is conformal to $\Bbb C$ is isometric to either the
%plane or a helicoid.\rm

\end{quote}

\section{Introduction}

A conformal orientation-preserving local diffeomorphism that is
defined in the open unit disc $D=\{z\in\Bbb C : |z|<1\}$ and takes
values into $\Bbb R^2$ can be viewed as a holomorphic function $f
: D\to\Bbb C$. Of special interest is the case when $f$ is
univalent, that is, injective. The class $S$ of all holomorphic
univalent functions
 in $D$ satisfying $f(0)=0$ and $f'(0)=1$ was the object of much study in the last century, culminating
 with the solution by de Branges (\cite{dB}, \cite{S}) of the celebrated Bieberbach conjecture: for any $f\in S$, the estimate
\begin{eqnarray}\label{Bie}
|f^{(k)}(0)|\leq k k!
\end{eqnarray}
holds for all $k\geq 2$. Equivalently, $|f^{(k)}(0)|\leq k k!
|f'(0)|$ for any injective holomorphic function on $D$. The case
$k=2$, due to Bieberbach, yields the so-called distortion theorems
which, in turn, imply the compactness of the class $S$ (\cite{Po},
\cite{S}). Thus, the basic estimate $|f''(0)|\leq 4$ for $f\in S$,
most commonly written in the form $|a_2|\leq 2$ where
$f(z)=z+a_2z^2+\cdots,$ already yields important qualitative
information. In particular, it follows from the compactness of $S$
that there are constants $C_k$ such that $|f^{(k)}(0)|\leq C_k
|f'(0)|$ for every $k\geq 2$ and injective holomorphic function
$f$ on $D$. The Bieberbach conjecture (the de Branges theorem)
asserts that one can take $C_k$ to be $k k!$.

The aim of this paper is to establish a  generalization of  Bieberbach's fundamental
estimate $|f''(0)|\leq 4 |f'(0)|$ that applies to  all  injective
 smooth conformal immersions
 $f : D\to\Bbb R^n,\,n\geq 2$.  The new
estimate involves two  correction terms. The
first one is geometric, coming from the second fundamental form of
the image surface $f(D)$. The second term
is of a dynamical nature,  and involves
certain Riemannian quantities associated  to conformal attractors.

Given a submanifold  of $\Bbb R^n$, we denote by $\nabla$ and $\sigma$  its  connection
and  second fundamental form, both extended by complex
linearity.  Recall also that a vector field on a Riemannian manifold is said to be
conformal if it generates a local flow of conformal maps. We use
subscripts to denote differentiation with respect to $z$, where
$2\,\partial_ z=\partial_ x-i\partial_ y$.

\vskip15pt

\begin{thm} Let $f: (D,0)\to (N,p)$ be an injective smooth conformal immersion,
where $N\subset\Bbb R^n$ is a smooth embedded topological disc,
$n\geq 3$. Let $\mathfrak X$ be the (non-empty) family of all
normalized conformal attractors on $(N,p)$, namely those smooth
vector fields $X$ on $N$ satisfying \vskip5pt \noindent{i)} X is
conformal,\;$X(p)=0$,\;$(\nabla X)_p=-I$.

\noindent{ii)} Every positive orbit of X tends to $p$.
\vskip5pt
\noindent{Then} {\em
\begin{eqnarray}\label{extbie}
\sup_{X\in\mathfrak
X}\big\Vert f_{zz}(0)-\sigma\big(f_z(0), f_z(0)\big)+(\nabla^2X)_p\big(f_z(0),f_z(0)\big)
\big\Vert\leq
4\Vert f_z(0)\Vert.\end{eqnarray}
\em}\end{thm}

\vskip10pt We expect  to develop further the ideas in this paper
and establish  analogues of the higher order estimates
$|f^{(k)}(0)|\leq kk! |f'(0)|$, $k\geq 3$,  in the broader context
of  injective  conformal immersions $f: D\to\mathbb R^n$. One
would then have, in all dimensions, a geometric-conformal version
of the de Branges theorem.

\vskip10pt
\noindent{\bf Example.} Taking $N=\mathbb R^2\subset\mathbb R^n$
as a totally geodesic submanifold and $X(w)=-(w-p)$,  we have $\nabla X=-I$  and so $\nabla^2 X=0$.
Since  $\sigma=0$,  one sees that the
original Bieberbach estimate $|f''(0)|\leq 4|f'(0)|$ can be
recovered from (\ref{extbie}).

\vskip10pt

It follows from  classical results that for any injective smooth
immersion $g:D\to \Bbb R^n$ there is a diffeomorphism $h:\Omega\to
D$, where $\Omega$ is either $D$ or $\Bbb C$, such that $f=g\circ
h$ is an injective conformal immersion. Restricting $h$ to $D$ in
the case $\Omega=\Bbb C$, one then sees  that  the above theorem
applies to every injective immersion, after a suitable
reparametrization.

We observe that the original Bieberbach estimate does not carry
over to the higher dimensional case. Indeed, if every injective
conformal immersion $f : D\to\mathbb R^n$ were to satisfy $\Vert
f_{zz}(0)\Vert\leq 4\Vert f_{z}(0)\Vert$, a contradiction could be
reached as follows. Let $g : \mathbb C\to\mathbb R^n$ be an
injective conformal harmonic immersion which is not totally
geodesic, i.e., $g(\mathbb C)$ is a parabolic simply-connected
embedded minimal surface. A concrete example is provided by the
helicoid given in coordinates $z=x+i y$ by  $g(x,y)=(\text{sinh}
x\, \text{cos} y, \text{sinh} x \,\text{sin} y, y)$ (see the
Conjecture below). Applying the above estimate to $f : D\to\mathbb
R^n,\, f(z)=g(Rz)$, and letting $R\to\infty$, one obtains
$$g_{zz}(0)=\frac{1}{4}\Big(g_{xx}(0)-2ig_{xy}(0)-g_{yy}(0)\Big)=\frac{1}{2}\Big(g_{xx}(0)-ig_{xy}(0)\Big)=0.$$ Replacing $g(z)$ by
$g(z+z_o)$ in the above reasoning, we obtain $g_{zz}\equiv 0$, and so $g_{xx}$, $g_{yy}$ and $g_{xy}$ are also identically zero,
contradicting the fact that $g(\mathbb C)$ is not a plane.

It is tempting to believe that the contribution in (\ref{extbie}) coming from the dynamic term can be made to vanish, regardless of the embedding,  in which case one would obtain  an
estimate that depends solely on the geometry of the surface
$f(D)$.  To investigate this possibility, suppose that $$\inf_{X\in\mathfrak
X}\big\Vert(\nabla^2X)_p\big(f_z(0),f_z(0)\big)\big\Vert=0$$ for
all conformal embeddings $f : D\to\mathbb R^n$ of the unit disc.
Consider the map $g$ above, parametrizing a  helicoid.  After
applying (\ref{extbie}) to $f : D\to\mathbb R^3,\,f(z)=g(Rz)$, and
reasoning as above, one has $$g_{zz}-\sigma\big(g_z,
g_z\big)=\frac{1}{2}\Big(g_{xx}^T-ig_{xy}^T\Big)\equiv 0.$$ In
particular, the coordinates curves of $g$ are geodesics of $g(D)$.
Since these curves are easily seen to be asymptotic lines of
$g(D)$, one would conclude that the traces of the coordinates
curves of $g$ are (segments of) straight lines in $\mathbb R^3$.
But this contradicts the fact that the coordinates curves
$y\mapsto g(x,y)$ are helices. Hence, the dynamic term is
essential for the validity of (\ref{extbie}).

As the reader may  have suspected by now,  part of our motivation
for proving a version of the Bieberbach estimate in the realm of
conformal embeddings $D\to \Bbb R^n$ comes  from the theory of
minimal surfaces.   We explain below how these two sets of ideas
merge.

The past few years have witnessed  great advances in the study of
simply-connected embedded minimal surfaces (i.e., minimal surfaces
without self-intersetions).  From an analytic standpoint, these
geometric objects correspond to conformal harmonic embeddings of
either $D$ or $\Bbb  C$ into $\Bbb R^3$.  In a series of
groundbreaking papers,  Colding and Minicozzi
(\cite{CM1}-\cite{CM5}) were able to give a very detailed
description  of the  structure of embedded minimal discs. Using
their theory, as well as other tools, Meeks and Rosenberg showed
in a landmark paper \cite{MR} that  helicoids and planes are the
only properly embedded simply-connected minimal surfaces in
$\mathbb R^3$ (subsequently, properness was weakened to mere
completeness).

The classical link between minimal surfaces and complex analysis
has been explored, with great success, to tackle other fundamental
geometric problems. Given the history of the subject, one is
naturally inclined to look for  a complex-analytic interpretation
of the works of Colding-Minicozzi and Meeks-Rosenberg, with the
hope that more could be revealed about the structure of embedded
minimal discs.  Although this effort is still in its infancy, one
can already  delineate the contours of a general programme. A
central theme to be explored is the role of the conformal type in
the embeddedness question for minimal surfaces. In particular, one
would like to know, in the Meeks-Rosenberg theorem, if
parabolicity alone suffices:

\vskip10pt

\noindent  \bf Conjecture. \rm    If $g:\Bbb C\to \Bbb R^3$ is a
conformal harmonic embedding, then $g(\mathbb C)$ is either a flat
plane or a    helicoid\rm.

\vskip10pt

There is a compelling analogy between the theory of conformal
harmonic  \it embeddings \rm of the open unit disc $ D \subset
\mathbb C $ into $\mathbb R^3$, and  the very rich theory of
holomorphic \it univalent \rm functions  on $D$.  It is an easy
matter to use a scaling argument, as it was done above, together
with the (classical) Bieberbach estimate to establish the scarcity
of  univalent entire functions: they are all of the form
$f(z)=az+b$, $ a\neq 0$. One ought to regard this statement as the
complex-analytic  analogue of the above conjecture.

A major hurdle in trying to use similar scaling  arguments to
settle the above conjecture is that one does not  have an \it a
priori \rm control on the dynamic term in  (\ref{extbie}).
Nevertheless, a first step  in the programme of using complex
analysis towards studying embedded minimal discs has been  taken
in \cite{FX}, where   a new proof was given of the classical
theorem of Catalan,  characterizing (pieces of) planes and
helicoids as the only ruled minimal surfaces.   The proof  of
Catalan's theorem is reduced, after careful normalizations,  to
the uniqueness of solutions of  certain    holomorphic differential
equations.

The present work is in fact part of a much larger programme, going
well beyond minimal surface theory, whose aim is to identify the
 analytic, geometric and topological mechanisms behind the
phenomenon of global injectivity; see, for instance,
\cite{AL}-\cite{BCW}, \cite{E}, \cite{G}, \cite{NTX}-\cite{Pl},
\cite{R}, \cite{X1}. The reader can find in \cite{X2}, p.17, a
short description of some of the recent results in the area of
global injectivity.

In closing, we would like to point out that the main result in
\cite{X1} --  a rigidity theorem characterizing the identity map
of $\mathbb C^n$ among injective local biholomorphisms --,  is
also  based on ideas suggested by the original Bieberbach
estimate.

\vskip20pt

\noindent {\bf Acknowledgements.} This paper  was written during an extended
visit of the first named author to the University of Notre Dame.
He would like to record his gratitude to the mathematics
department for the invitation,
as well as for its hospitality.

\section{Proof of the generalized Bieberbach estimate}
\vskip5pt
The proof of Theorem 1.1 will be split into a series of lemmas,
reflecting the dynamical, complex-analytic and Riemannian aspects
of the argument.
\vskip10pt
\begin{lem}\label{dyn} Let $\widetilde X$ be a smooth
vector field defined in an open set $U\subset \Bbb R^n$. If
$\widetilde X(p)=0$ and $(d\widetilde X)_p= -I$, then the local
flow $\eta_t$ of $\widetilde X$ satisfies

\vskip10pt

\noindent{\em (i)} $(d\eta_t)_p=e^{-t} I,\;t>0$.

\noindent{\em (ii)} $\displaystyle \lim_{t\to\infty}\Vert
(d\eta_t)_p\Vert^{-1} (d^2\eta_t)_p (v, w)=(d^2\widetilde X)_p (v,
w),\;\;v, w\in\mathbb R^n.$
\end{lem}
\vskip10pt
\noindent{\bf Proof. (i)} For all $x\in U$ and $v\in \Bbb R^n$,
\begin{eqnarray}\label{din1}
\frac{d}{dt}(d\eta_t)_x(v)&=&\frac{d}{dt}\left\{\frac{d}{ds}\Big|_{s=0}\eta_t(x+sv)\right\}
=\frac{d}{ds}\Big|_{s=0}\left(\frac{d}{dt}\eta_t(x+sv)\right)\nonumber\\
&=&\frac{d}{ds}\Big|_{s=0}\widetilde X\circ\eta_t
(x+sv)=(d\widetilde X)_{\eta_t(x)}\circ (d\eta_t)_x (v).
\end{eqnarray}
Since $\widetilde X(p)=0$, $\eta_t(p)=p$, for all t. Setting $x=p$
in (\ref{din1}), one obtains
\begin{eqnarray}\label{din2}
\frac{d}{dt}(d\eta_t)_p(v)=(d\widetilde X)_p\circ (d\eta_t)_p
(v)=-(d\eta_t)_p(v).
\end{eqnarray}
Since $\eta_0(x)=x$ for all $x$, $(d\eta_0)_p=I$ and, by
(\ref{din2}),
$$
(d\eta_t)_p(v)=e^{-t}(d\eta_0)_p(v)=e^{-t}v,\;\;v\in\Bbb R^n.
$$

\noindent{\bf (ii)} For all $x\in U$ and $v, w\in \Bbb R^n$,
\begin{eqnarray}
(d^2\eta_t)_x(v, w)= \frac{d}{ds}\Big|_{s=0}\left\{
\frac{d}{dh}\Big|_{h=0} \eta_t(x+sv+hw)\right\}.
\end{eqnarray}
Thus,
\begin{eqnarray}\label{din3}
\frac{d}{dt}(d^2\eta_t)_x(v, w)&=& \frac{d}{ds}\Big|_{s=0}\left\{
\frac{d}{dh}\Big|_{h=0}\,
\frac{d}{dt}\eta_t(x+sv+hw)\right\}\nonumber\\&=&
\frac{d}{ds}\Big|_{s=0}\left\{\frac{d}{dh}\Big|_{h=0}\widetilde
X\circ\eta_t(x+sv+hw)\right\}\nonumber\\&=&d^2(\widetilde
X\circ\eta_t)_x(v, w).
\end{eqnarray}
It follows that
\begin{eqnarray}\label{din4}
\frac{d}{dt}(d^2\eta_t)_x(v, w)=(d^2\widetilde
X)_{\eta_t(x)}\big((d\eta_t)_x v, (d\eta_t)_x w\big)+(d\widetilde
X)_{\eta_t(x)}\big((d^2\eta_t)_x(v, w)\big).
\end{eqnarray}
Taking $x=p$ in (\ref{din4}) and recalling that $(d\widetilde
X)_p=-I$, we have
\begin{eqnarray}\label{din5}
\frac{d}{dt}(d^2\eta_t)_p(v, w)&=&(d^2\widetilde
X)_p\big((d\eta_t)_p v, (d\eta_t)_p w\big)+(d\widetilde
X)_p\big((d^2\eta_t)_p(v, w)\big)\nonumber\\&=&(d^2\widetilde
X)_p\big((d\eta_t)_p v, (d\eta_t)_p w\big)-(d^2\eta_t)_p(v, w).
\end{eqnarray}
From (i) and (\ref{din5}), one has
\begin{eqnarray}
\frac{d}{dt}(d^2\eta_t)_p(v, w)=e^{-2 t}(d^2\tilde X)_p(v,
w)-(d^2\eta_t)_p(v, w).
\end{eqnarray}
Letting $V(t)=(d^2\eta_t)_p(v, w)$ and $W=(d^2\tilde X)_p(v, w)$,
the last equation becomes
\begin{eqnarray}\label{din6}
\frac{d}{dt}V(t)=-V(t)+e^{-2 t}W,
\end{eqnarray}
so that
\begin{eqnarray}
\frac{d}{dt}(e^tV(t))=e^{-t}W,
\end{eqnarray}
and, by integration,
\begin{eqnarray}\label{din7}
e^tV(t)-V(0)=W(1-e^{-t}).
\end{eqnarray}
Since $\eta_0(x)=x,\,x\in U$, it follows that
$V(0)=(d^2\eta_0)_p(v, w)=0$. Using this in (\ref{din7}), we have
\begin{eqnarray}\label{din8}
e^t V(t)=W(1-e^{-t}).
\end{eqnarray}
Taking the limit as $t\to\infty$ and using (i) one obtains (ii).
\qed \vskip20pt

\begin{lem} Let $f : D\rightarrow \Bbb R^n$ and $g : D\rightarrow \Bbb
R^n$ be two conformal embeddings of the unit disc into $\Bbb R^n$
such that $g(D)\subset f(D)$ and $f(0)=g(0)$. Assume that the
orientations induced on $g(D)$ by $f$ and $g$ are the same. Then there exists $\zeta\in\mathbb C,\;|\zeta|=1$, such that
\begin{eqnarray}\label{din9}
\left\Vert\frac{g_{zz}(0)}{\Vert
g_z(0)\Vert^2}-\zeta\frac{f_{zz}(0)}{\Vert
f_z(0)\Vert^2}\right\Vert\leq \frac{4}{\Vert g_z(0)\Vert}.
\end{eqnarray}

\end{lem}
\vskip10pt
\noindent{\bf Proof.} The function $\varphi=f^{-1}\circ g: D\to D$
is well defined, conformal and orientation preserving. Hence
$\varphi$ is holomorphic, $\varphi (0)=0$. Using $g_z(0)=f_z(0)
\varphi_z(0)$ and
$g_{zz}(0)=f_{zz}(0)\varphi_z(0)^2+f_z(0)\varphi_{zz}(0)$, one
computes
\begin{eqnarray}
\frac{g_{zz}(0)}{\Vert g_z(0)\Vert^2}=\frac{f_{zz}(0)}{\Vert
f_z(0)\Vert^2}\,\frac{\varphi_z(0)^2}{|\varphi_z(0)|^2}+\frac{f_z(0)\,
\varphi_{zz}(0)}{\Vert f_z(0)\Vert^2\,|\varphi_z(0)|^2}.
\end{eqnarray}
But
\begin{eqnarray}
\left\Vert \frac{f_z(0)\, \varphi_{zz}(0)}{\Vert
f_z(0)\Vert^2\,|\varphi_z(0)|^2}\right\Vert\leq\frac{4}{\Vert
f_z(0)\Vert |\varphi_z(0)|}=\frac{4}{\Vert g_z(0)\Vert},
\end{eqnarray}

\noindent{by} Bieberbach's inequality $|\psi''(0)|\leq 4
|\psi'(0)|$, valid for all univalent holomorphic functions $\psi:
D\to \Bbb C$. Hence
\begin{eqnarray}\label{fd}
\left\Vert\frac{g_{zz}(0)}{\Vert
g_z(0)\Vert^2}-\frac{\varphi_z(0)^2}{|\varphi_z(0)|^2}\,\frac{f_{zz}(0)}{\Vert
f_z(0)\Vert^2}\right\Vert\leq\frac{4}{\Vert g_z(0)\Vert}
\end{eqnarray}

\noindent{and the} lemma follows from (\ref{fd}) with
$\zeta=\frac{\varphi_z(0)^2}{|\varphi_z(0)|^2}$.\qed

\newpage

\begin{lem} Let $N\subset \Bbb R^n,\,n\geq 2$, be a smooth surface and $f : (D, 0)\to (N,p)$
an injective smooth conformal immersion. Let $U$ be an open neighborhood
of $N$ in $\Bbb R^n$ and $\widetilde X$ a smooth vector field on
$U$ such that:

\vskip5pt

\noindent{\em (i)} $\widetilde X(p)=0$, $(d\widetilde X)_p=-I$.

\noindent{\em (ii)} $\widetilde X$ is tangent to $N$, $\widetilde
X|_N$ is conformal and every positive orbit of $\widetilde X|_N$
converges to $p$.

\vskip5pt

\noindent{Then}
\begin{eqnarray}\label{des}
\big\Vert(d^2\widetilde X)_p\big(f_z(0),
f_z(0)\big)+f_{zz}(0)\big\Vert\leq4\Vert f_z(0)\Vert.
\end{eqnarray}
\end{lem}

\vskip10pt

\noindent{\bf Proof.} We may assume  that $f(D)$ has compact closure in $N$ (the general case follows by replacing $f(z)$ by $f(Rz)$, with $R<1$, and letting $R\to 1$). Since every positive orbit of $\widetilde X|_N$
converges to $p$, there exists $t_0>0$ such that  the positive flow $\eta_t$ of $\widetilde X$ satisfies $\eta_t(f(D))\subset f(D)$ for all $t\geq t_0$.
 Since $\widetilde X|_N$ is conformal, so is
$g^{(t)}=g=\eta_t\circ f$.   It follows easily from Lemma 2.1 (i) that  $f$ and $g$
induce the same orientations on $g(D)\subset f(D)$, as required by Lemma 2.2. We have
\begin{eqnarray}\label{din10}
g_z(z)=(\eta_t\circ f)_z(z)=(d\eta_t)_{f(z)}(f_z(z))
\end{eqnarray}
and
\begin{eqnarray}\label{din11}
g_{zz}(0)=(d^2\eta_t)_p(f_z(0), f_z(0))+(d\eta_t)_p(f_{zz}(0)).
\end{eqnarray}

\noindent{By} (\ref{din9}),
\begin{eqnarray}\label{din12}
\frac{\Vert g_{zz}(0)\Vert}{\Vert g_z(0)\Vert^2}\leq
\frac{4}{\Vert g_z(0)\Vert}+\frac{\Vert f_{zz}(0)\Vert}{\Vert
f_z(0)\Vert^2}.
\end{eqnarray}
From (\ref{din10}), (\ref{din11}) and (\ref{din12}), we then have
\begin{eqnarray}
 \frac{\Vert (d^2\eta_t)_p(f_z(0),
f_z(0))+(d\eta_t)_p(f_{zz}(0))\Vert}{\Vert
(d\eta_t)_p\Vert^2\,\Vert f_z(0)\Vert^2}\leq\frac{4}{\Vert
(d\eta_t)_p\Vert\,\Vert f_z(0)\Vert}+\frac{\Vert
f_{zz}(0)\Vert}{\Vert f_z(0)\Vert^2}
\end{eqnarray}

\noindent{Multiplying} the last equation by $\Vert
(d\eta_t)_p\Vert\Vert f_z(0)\Vert^2$,
\begin{eqnarray}
 \frac{\Vert (d^2\eta_t)_p(f_z(0),
f_z(0))+(d\eta_t)_p(f_{zz}(0))\Vert}{\Vert (d\eta_t)_p\Vert}\leq 4
\Vert f_z(0)\Vert+\Vert (d\eta_t)_p\Vert\,\Vert f_{zz}(0)\Vert.
\end{eqnarray}

\noindent{Taking} the limit as $t\to\infty$ and using Lemma
\ref{dyn}, one obtains (\ref{des}).\qed

\vskip20pt

In the next lemma we denote by $\overline\nabla$ and $\nabla$ the
Riemannian connections of $\mathbb R^n$ and $M$, respectively, and
by $\sigma$ the second fundamental form of $M$.

\vskip15pt

\begin{lem}Let $M\subset\Bbb R^n$ be an $m$-dimensional  submanifold and $X$  a
smooth vector field on $M$ such that $X(p)=0$ and $\nabla_v X=-v$,
for some $p\in M$ and all $v\in T_pM$. Then, for any extension
$\widetilde X$ of $X$ to an open neighborhood of $p$ in $\Bbb R^n$
and all $v, w\in T_pM$, one has
\begin{eqnarray}\label{int1}
(d^2\widetilde
X)_p(v,w)=(\nabla^2X)_p(v,w)-\overline{\nabla}_{\sigma (v,
w)}\,\widetilde X-2 \sigma (v, w).
\end{eqnarray}
\end{lem}
\vskip10pt
\noindent{\bf Proof.} Writing $u_1,\dots,u_n$ for the canonical
basis of $\mathbb R^n$  one has,  for all $v, w\in \Bbb R^n$,
\begin{eqnarray}\label{int2}
(d^2\widetilde X)_p(v, w)=\sum_{i,
j=1}^n\frac{\partial^2\widetilde X}{\partial x_i\partial
x_j}(p)\langle v,u_i\rangle\langle w,u_j\rangle,
\end{eqnarray}
so that the coordinates of this vector are given by
\begin{eqnarray}\label{int3}
\sum_{i, j=1}^n\frac{\partial^2\widetilde X_k}{\partial
x_i\partial x_j}(p)\langle v,u_i\rangle\langle
w,u_j\rangle&=&\sum_{i, j=1}^n\langle v,u_i\rangle\langle
w,u_j\rangle\text{hess} \widetilde X_k(p) (u_i,
u_j)\nonumber\\&=&\text{hess}\widetilde X_k(p)(v,
w),\;\;k=1,\dots,n.
\end{eqnarray}
Choose an orthonormal frame field $\{e_1,\dots,e_n\}$ in an open
neighborhood of $p$ in such a way that $e_1,\dots,e_m$ are tangent
along $M$, $e_{m+1},\dots,e_n$ are normal along $M$ and
$\nabla_{e_i}e_j(p)=0,\,i, j=1,\dots,m$. Hence, by (\ref{int3}),

\begin{eqnarray}\label{int4}
\sum_{i, j=1}^n\frac{\partial^2\widetilde X_k}{\partial
x_i\partial x_j}(p)\langle v,u_i\rangle\langle
w,u_j\rangle&=&\sum_{i, j=1}^nv_iw_j\text{hess} \widetilde X_k(p)
(e_i, e_j)\nonumber\\&=&\sum_{i, j=1}^n\langle
\overline{\nabla}_{e_i}\text{grad} \widetilde X_k (p),
e_j(p)\rangle v_i w_j,
\end{eqnarray}
where the $v_i's$ and $w_i's$ are the components of $v$ and $w$ in
the basis $\{e_1,\dots,e_n\}$. It follows that
\begin{eqnarray}\label{int5}
\sum_{i, j=1}^n\frac{\partial^2\widetilde X_k}{\partial
x_i\partial x_j}(p)\langle v,u_i\rangle\langle
w,u_j\rangle&=&\sum_{i, j=1}^nv_i w_j\big(e_i\langle \text{grad}
\widetilde X_k , e_j\rangle-\langle \text{grad} \widetilde X_k ,
\overline{\nabla}_{e_i}e_j\rangle\big)(p) \nonumber\\&=&\sum_{i,
j=1}^nv_i w_j\big(e_i e_j (\widetilde X_k)-\overline{\nabla}_{e_i}
e_j (\widetilde X_k)\big)(p).
\end{eqnarray}
From (\ref{int2}) and (\ref{int5}) we obtain,
\begin{eqnarray}\label{int6}
(d^2\widetilde X)_p(v, w)=\sum_{i,
j=1}^n\big(\overline{\nabla}_{e_i}\overline{\nabla}_{e_j}\widetilde
X\big)(p) v_i w_j-\sum_{i,
j=1}^n\big(\overline{\nabla}_{\overline{\nabla}_{e_i}
e_j}\widetilde X\big)(p) v_i w_j.
\end{eqnarray}
Suppose now that $v$ and $w$ are tangent to $M$. Since the
restriction to $M$ of $\{e_1,\dots,e_n\}$ is an adapted frame
satisfying $\nabla_{e_i}e_j (p)=0,\,i, j=1,\dots,m$, and
$\widetilde X$ restricted to $M$ is $X$, it follows from
(\ref{int6}) and the Gauss equation \cite{Dj}
\begin{eqnarray}\label{GE}
\overline\nabla_YX=\nabla_YX+\sigma(X,Y)
\end{eqnarray}
that
\begin{eqnarray}\label{int7}
(d^2\widetilde X)_p(v, w)&=&\sum_{i,
j=1}^m\big(\overline{\nabla}_{e_i}\overline{\nabla}_{e_j}
X\big)(p) v_i w_j-\sum_{i, j=1}^m\big(\overline{\nabla}_{\sigma
(e_i, e_j)}\widetilde X\big) v_i w_j\nonumber\\&=&\sum_{i,
j=1}^m\big(\overline{\nabla}_{e_i}\overline{\nabla}_{e_j}
X\big)(p) v_i w_j-\overline{\nabla}_{\sigma(v, w)}\widetilde X.
\end{eqnarray}
We will now determine the first term on the right hand side of
(\ref{int7}). From
\begin{eqnarray}
\overline{\nabla}_{e_j}X=\nabla_{e_j}X+\sigma(e_j, X),\nonumber
\end{eqnarray}
we have at $p$, for $1\leq i, j\leq m$,
\begin{eqnarray}
\overline{\nabla}_{e_i}\overline{\nabla}_{e_j}X=\overline{\nabla}_{e_i}\nabla_{e_j}X+\overline{\nabla}_{e_i}\sigma(e_j,
X)=\nabla_{e_i}\nabla_{e_j}X+\sigma(e_i,
\nabla_{e_j}X)+\overline{\nabla}_{e_i}\sigma(e_j, X).\nonumber
\end{eqnarray}
Hence
\begin{eqnarray}\label{int8}
\sum_{i,
j=1}^m(\overline{\nabla}_{e_i}\overline{\nabla}_{e_j}X)(p)v_i
w_j&=&\sum_{i, j=1}^m(\nabla_{e_i}\nabla_{e_j}X)(p)v_i
w_j+\sum_{i, j=1}^m\sigma(e_i, \nabla_{e_j}X)(p)v_i
w_j\nonumber\\&+&\sum_{i,
j=1}^m\big(\overline{\nabla}_{e_i}\sigma(e_j, X)\big)(p)v_i w_j.
\end{eqnarray}

\noindent{The} third term on the right hand side of (\ref{int8})
involves the quantity $\overline{\nabla}_{e_i}\sigma(e_j, X)$
which can be computed using the normal connection
$\nabla^{\perp}$:
\begin{eqnarray}\label{int9}
\overline{\nabla}_{e_i}\sigma(e_j,
X)&=&\nabla_{e_i}^{\bot}\sigma(e_j,
X)+(\overline{\nabla}_{e_i}\sigma(e_j,
X))^T\nonumber\\&=&(\nabla_{e_i}^{\bot}\sigma)(e_j, X)+\sigma(e_j,
\nabla_{e_i}X)+\sum_{k=1}^m\big\langle
(\overline{\nabla}_{e_i}\sigma(e_j, X))^T, e_k\big\rangle
e_k\nonumber\\&=&(\nabla_{e_i}^{\bot}\sigma)(e_j, X)+\sigma(e_j,
\nabla_{e_i}X)+\sum_{k=1}^m\big\langle
\overline{\nabla}_{e_i}\sigma(e_j, X), e_k\big\rangle
e_k\nonumber\\&=&(\nabla_{e_i}^{\bot}\sigma)(e_j, X)+\sigma(e_j,
\nabla_{e_i}X)-\sum_{k=1}^m\big\langle\sigma(e_j, X),
\overline\nabla_{e_i}e_k\big\rangle
e_k\nonumber\\&=&(\nabla_{e_i}^{\bot}\sigma)(e_j, X)+\sigma(e_j,
\nabla_{e_i}X)-\sum_{k=1}^m\big\langle\sigma(e_j, X), \sigma(e_i,
e_k)\big\rangle e_k.
\end{eqnarray}
From (\ref{int8}) and (\ref{int9}),
\begin{eqnarray}\label{int10}
\sum_{i,
j=1}^m(\overline{\nabla}_{e_i}\overline{\nabla}_{e_j}X)v_i
w_j&=&\sum_{i, j=1}^m(\nabla_{e_i}\nabla_{e_j}X)v_i w_j+\sigma(v,
\nabla_wX)+\sum_{i, j=1}^m(\nabla_{e_i}^{\bot}\sigma)(e_j, X)v_i
w_j\nonumber\\&+&\sum_{i, j=1}^m\sigma(e_j, \nabla_{e_i}X)v_i
w_j-\sum_{i, j, k=1}^m\Big(\big\langle\sigma(e_j, X), \sigma(e_i,
e_k)\big\rangle e_k\Big)v_i w_j\nonumber\\&=&\sum_{i,
j=1}^m(\nabla_{e_i}\nabla_{e_j}X)v_i w_j+\sigma(v,
\nabla_wX)+(\nabla_v^{\bot}\sigma)(w, X)\nonumber\\&+&\sigma(w,
\nabla_vX)-\sum_{k=1}^m\langle \sigma(w,X),\sigma(v,e_k)\rangle
e_k.
\end{eqnarray}

\noindent{Since}, by assumption, $X(p)=0$ and $(\nabla X)_p=-I$,
we have
\begin{eqnarray}\label{int11}
\sum_{i,
j=1}^m(\overline{\nabla}_{e_i}\overline{\nabla}_{e_j}X)(p)v_i
w_j=\sum_{i, j=1}^m(\nabla_{e_i}\nabla_{e_j}X)(p)v_i
w_j-2\sigma(v, w).
\end{eqnarray}
It follows from (\ref{int7}) and (\ref{int11}) that
\begin{eqnarray}\label{int12}
(d^2\widetilde X)_p(v, w)=\sum_{j=1}^m(\nabla_v\nabla_{e_j}X)(p)
w_j-2\sigma(v, w)-\overline{\nabla}_{\sigma(v, w)}\widetilde X.
\end{eqnarray}
Let $V$ and $W$ be arbitrary smooth extensions of $v$ and $w$ to
an open neighborhood of $p$ in $M$. Using again $X(p)=0$ and
$(\nabla X)_p=-I$ and recalling that $\nabla_{e_j}e_i(p)=0$, we
have, from the definition of the curvature tensor $R$,
\begin{eqnarray}\label{int13}
\sum_{j=1}^m(\nabla_v\nabla_{e_j}X)(p) w_j&=&\sum_{j=1}^m \big(R(V,
e_j)X+\nabla_{e_j}\nabla_VX+\nabla_{[V,
e_j]}X\big)(p)w_j\nonumber\\&=&\nabla_W\nabla_VX(p)+\nabla_WV(p).
\end{eqnarray}

\noindent{From} (\ref{int12}) and (\ref{int13}) we obtain
\begin{eqnarray}\label{int14}
(d^2\widetilde X)_p(v, w)=\nabla_W\nabla_V X(p)+\nabla_W V
(p)-\overline{\nabla}_{\sigma (v, w)}\,\widetilde X-2 \sigma (v,
w).
\end{eqnarray}

\noindent{On} the other hand,
\begin{eqnarray}\label{int15}
(\nabla^2X)_p (v,w)&=&\big(\nabla_w\nabla
X\big)(v)=\big(\nabla_W\nabla X(V)\big)(p)-\nabla
X\big(\nabla_WV\big)(p)\nonumber\\&=&\nabla_W\nabla_VX(p)-\nabla_{\nabla_WV}X(p).
\end{eqnarray}

\noindent{Since} $(\nabla X)_p=-I$, we thus obtain
\begin{eqnarray}\label{int16}
(\nabla^2X)_p (v,w)=\nabla_W\nabla_V X(p)+\nabla_WV(p).
\end{eqnarray}

\noindent{Formula} (\ref{int1}) now follows from (\ref{int12}),
(\ref{int13}) and (\ref{int16}).\qed

\vskip20pt

\noindent{\bf Proof of Theorem 1.1.} We begin by showing that
$\mathfrak X\neq\emptyset$. Being simply-connected, the Riemannian
surface $N$ is globally conformally flat and so it is isometric to
$(\Omega, \widetilde g)$, where $\Omega$ is either $D$ or $\mathbb
C$, $\widetilde g=e^{2 \varphi}g$ for some smooth function
$\varphi$, and $g$ is the standard flat metric on $\Omega$.
Writing $\widetilde\nabla$ and $\nabla$ for the Riemannian
connections of $\widetilde g$ and $g$, respectively, one has
(\cite{Pr} p. 172):
\begin{eqnarray}\label{cm1}
\widetilde\nabla_YX=\nabla_YX+Y(\varphi)X+X(\varphi)Y-\langle
X,Y\rangle\nabla\varphi.
\end{eqnarray}
Being holomorphic, $X(z)=-z$ generates a flow of conformal maps
that clearly leaves $\Omega$ invariant. It now follows from
(\ref{cm1}) that $\widetilde\nabla_YX=\nabla_YX=-Y$ at $0$, thus
showing that $X\in \mathfrak X$.

We now proceed to establish the estimate (\ref{extbie}). Since $N$
is contractible, one can choose a global orthonormal frame
$\{\xi_1,\dots,\xi_{n-2}\}$ of the normal bundle of $N$ in
$\mathbb R^n$. Extend $X$ to an open neighborhood of $N$ in $\Bbb
R^n$ by
\begin{eqnarray}\label{int17}
\widetilde
X\left(q+\sum_{i=1}^{n-2}s_i\xi_i(q)\right)=X(q)-\sum_{i=1}^{n-2}s_i\xi_i(q),\;\;
\;|s_i|<\epsilon (q),
\end{eqnarray}

\noindent{with} $\epsilon(q)$ sufficiently small, $q\in N$. For
all $v\in T_pN$ we have
$$
(d\widetilde
X)_pv=\overline\nabla_v\widetilde X=\overline\nabla_vX=\nabla_vX+\sigma(v,
X(p))=-v,
$$
since $X(p)=0$ and $(\nabla X)_p=-I$. On the other hand,
\begin{eqnarray}
\overline\nabla_{\xi_i(p)}\widetilde X=(d\widetilde
X)_p\xi_i(p)&=&\frac{d}{ds}\Big|_{s=0}\widetilde X\big(p+s
\xi_i(p)\big)\nonumber\\&=&\frac{d}{ds}\Big|_{s=0}\big(X(p)-s\xi_i(p)\big)=-\xi_i(p).\nonumber
\end{eqnarray}
It follows from the last two equations that
$(\overline\nabla\widetilde X)_p=(d\widetilde X)_p=-I$. By Lemma
2.3, we have
\begin{eqnarray}\label{int18}
\big\Vert(d^2\widetilde X)_p\big(f_z(0),
f_z(0)\big)+f_{zz}(0)\big\Vert\leq 4\Vert f_z(0)\Vert.
\end{eqnarray}
Also, from Lemma 2.4 and $(\overline{\nabla}\widetilde X)_p=-I$,
\begin{eqnarray}\label{int19}
(d^2\widetilde X)_p(v,
w)&=&(\nabla^2X)_p(v,w)-\overline{\nabla}_{\sigma (v,
w)}\,\widetilde X-2 \sigma (v, w),
\nonumber\\&=&(\nabla^2X)_p(v,w)- \sigma (v, w),\;\; v, w\in T_pN,
\end{eqnarray}
which, in particular, implies that $(\nabla^2 X)_p$ is symmetric.
From $f_z(0)=\frac{1}{2}\big(f_x(0)-if_y(0)\big)$ and
\begin{eqnarray}\label{int20}
(d^2\widetilde X)_p\big(f_z(0),
f_z(0)\big)&=&\frac{1}{4}(d^2\widetilde X)_p\big(f_x(0)-if_y(0),
f_x(0)-if_y(0)\big) \nonumber\\&=&\frac{1}{4}\Big\{(d^2\widetilde
X)_p\big(f_x(0), f_x(0)\big)-(d^2\widetilde X)_p\big(f_y(0),
f_y(0)\big)\nonumber\\&-&2i (d^2\widetilde X)_p\big(f_x(0),
f_y(0)\big)\Big\},
\end{eqnarray}
we obtain, applying (\ref{int19}) to each term in (\ref{int20})
and rearranging the terms,
\begin{eqnarray}\label{int21}
(d^2\widetilde X)_p\big(f_z(0),
f_z(0)\big)=(\nabla^2X)_p\big(f_z(0),f_z(0)\big)- \sigma
\big(f_z(0), f_z(0)\big).
\end{eqnarray}

\noindent{Formula} (\ref{extbie}) follows from (\ref{int18}) and
(\ref{int21}).\qed

\noindent Francisco Fontenele\\Departamento de Geometria\\
Universidade Federal Fluminense\\Niter\'oi, RJ,
Brazil\\fontenele@mat.uff.br\\ \and \\
Frederico Xavier\\Department
of Mathematics\\University of Notre Dame\\Notre Dame, IN,
USA\\fxavier@nd.edu

\end{document}